\documentclass[fleqn,10 pt]{article}
 \usepackage[centertags]{amsmath}
\usepackage{amsfonts}
\usepackage{amssymb}
\usepackage{amsthm}
\usepackage{multirow}
\usepackage{newlfont}
\usepackage{graphicx}
\usepackage{epstopdf}
\usepackage{mathrsfs}
\usepackage[english]{babel}
 \usepackage{cite}
 \usepackage{color}
  \usepackage{makecell}
\newcommand{\q}{\quad}
\newcommand{\qq}{\qquad}

\newcommand{\bq }{\begin{equation}}
\newcommand{\eq }{\end{equation}}
\newcommand{\bbb }{\begin{eqnarray}}
\newcommand{\eee }{\end{eqnarray}}
\newcommand{\bb }{\begin{eqnarray*}}
\newcommand{\ee }{\end{eqnarray*}}
\newcommand{\ed }{\end{document}}
\theoremstyle{plain}
\newtheorem{thm}{Theorem}[section]

\newtheorem{rem}[thm]{Remark}
\newtheorem{example}{ Example}

\theoremstyle{definition}

\theoremstyle{example}

\numberwithin{equation}{section}

\begin{document}
\title{New recursive approximations for variable-order fractional operators with applications}

\author{M.A. Zaky$^{1}$, E.H. Doha$^{2}$, T.M. Taha$^{3}$,  D. Baleanu$^{4,5}$\\
	\footnotesize{$^{1}$Department of Applied Mathematics, National Research Centre
Dokki, 12622 Giza, Egypt}  \\
\footnotesize{$^{2}$Department of Mathematics, Faculty of Science, Cairo University
Giza, Egypt}  \\
\footnotesize{$^{3}$Department of Mathematics, Faculty of Science, Beni-Suef University
Beni-Suef, Egypt}  \\
\footnotesize{$^{4}$Department of Mathematics, Cankaya University
Ankara, Turkey}  \\
\footnotesize{$^{5}$nstitute of Space Sciences
Magurele-Bucharest, Romania}  \\
	\footnotesize{E-mail(corresp.):  ma.zaky@yahoo.com
	}}
	
\date{}
\maketitle

\begin{abstract}
 To broaden the range of applicability of variable-order fractional  differential models, reliable numerical approaches are needed to solve the model equation. In this paper, we develop Laguerre spectral collocation methods for solving variable-order fractional initial value problems on the half line.  Specifically, we derive three-term recurrence relations  to efficiently calculate the variable-order fractional integrals and derivatives of the modified generalized  Laguerre polynomials, which lead to the corresponding fractional
differentiation matrices that will be used to construct the collocation methods. Comparison with other existing methods shows the superior accuracy
of the proposed spectral collocation methods. \end{abstract}



$ $

\hrule

$ $

\section{Introduction}
The variable-order fractional (VO-F) operators \cite{samko1993integration,coimbra2003mechanics}, which are  generalizations  of constant-order
fractional operators \cite{zaky2017Legendre}, open up new possibilities for robust mathematical
modeling and simulation of diverse physical problems in science and engineering,  such as  modeling of
diffusive-convective effects on the oscillatory flows
\cite{pedro2008variable},
linear and nonlinear oscillators
with viscoelastic damping \cite{coimbra2003mechanics}, processing of geographical data using
VO-F derivatives \cite{cooper2004filtering},  constitutive laws in viscoelastic continuum
mechanics \cite{ramirez2007variable}, signature verification
through variable/adaptive fractional order differentiators
\cite{tseng2006design}, anomalous diffusion problems
\cite{zaky2016efficient,fu2015method} and chloride ions sub-diffusion in concrete
structures \cite{wei2017time}.  The VO-F operators can be employed to depict the
variable memory of systems \cite{lorenzo2002variable}.

The VO-F operators are nonlocal with singular kernels, which
makes the VO-F models complicated. Hence, the
solution of VO-F models is also more complicated. Numerical computation of the VO-F operators  is the key to understand
the behavior and physical meaning of the VO-F  models. Fu et al. \cite{fu2015method} applied
the method of approximate particular solutions to VO-F diffusion models. Cao and Qiu \cite{cao2016high} proposed a second order numerical approximation
via the VO-F weighted and shifted Gr\"{u}nwald-Letnikov formula to VO-F Riemann-Liouville derivative, and used it to solve VO-F ordinary differential equations. Zayernouri and Karniadakis \cite{zayernouri2015fractional} introduced fractional spectral collocation methods for
linear and nonlinear VO-F  differential equations.  Atangana et al. \cite{atangana2015stability} developed the Crank-Nicholson scheme to handle the time VO-F telegraph equation. Bhrawy et al. proposed accurate  spectral collocation methods for  VO-F differential equations such as  Schr\"{o}dinger equation \cite{bhrawy2017improved},  Galilei invariant advection diffusion equations \cite{abd2017space},   diffusion equation \cite{zaky2016efficient} and cable equation \cite{bhrawy2015numerical}. Moghaddam et al \cite{moghaddam2017extended,moghaddam2017integro}  developed accurate and robust algorithms for approximating VO-F derivatives and integrals. Tayebi et al \cite{tayebi2017meshless} proposed an accurate and robust meshless method based on the moving least squares approximation and the finite difference scheme for the numerical solution of VO-F advection-diffusion equation  on two-dimensional  arbitrary domains.

In this paper, we focus on the computation of the VO-F integrals and derivatives of the modified generalized Laguerre polynomials.
Applications of the constructed computations are illustrated to compute the VO-F Caputo   derivative. Besides, using the modified generalized Laguerre polynomials as the basis functions, we develop Laguerre-Gauss collocation  methods to solve fractional differential equations  of variable and constant orders on the half line.

This paper  is organized as  follows. Section \ref{sec2} presents the  fundamentals of  VO-F operators and properties of the modified generalized Laguerre polynomials. Numerical algorithms for
calculating the VO-F integral and the Caputo derivative are presented
in Sections \ref{sec3} and \ref{sec4}, respectively. The applications of the algorithms are
illustrated in Section \ref{sec5}. Numerical examples are presented in Section \ref{sec6}, and the conclusion
is drawn in the last section.

\section{Preliminaries and fundamentals} \label{sec2}

In this section, we concisely point out some definitions  of the VO-F operators \cite{coimbra2003mechanics,lorenzo2002variable,samko1993integration}. We  then  collect  some  important properties of the modified generalized Laguerre polynomials \cite{baleanu2013modified}. Assume that  $u(x)=0 $ for $x<0$.

\begin{enumerate}
  \item  The following VO-F integral operator was proposed in \cite{samko1993integration}
 \begin{equation}\label{za1}
 {}_0I_x^{\varrho  (x)} [u]: = x \mapsto \frac{1}
{{\Gamma (\varrho (x))}}\int_0^x {(x - r)^{\varrho (x) - 1} u(r)dr},\q  x \geq 0.\end{equation}
  \item In \cite{lorenzo2002variable} several definitions were proposed. The first is identical to \eqref{za1}.
The next one, is
\begin{equation}\label{za2}
{}_0I_x^{\varrho (x)} [u]: = x \mapsto \int_0^x {\frac{{(x - r)^{\varrho (r) - 1} }}
{{\Gamma (\varrho (r))}}u(r)dr},\q  x \geq 0.\end{equation}
  \item The following  operator was introduced in \cite{lorenzo2002variable}
\begin{equation}\label{za3}
{}_0I_x^{\varrho (x)} [u]: = x \mapsto \int_0^x {\frac{{(x - r)^{\varrho (x - r) - 1} }}
{{\Gamma (\varrho (x - r))}}u(r)dr}\q  x \geq 0.\end{equation}
\end{enumerate}

 The VO-F Caputo  derivative  could now be defined, as in the case of constant order \cite{samko1993integration}, as follows
\begin{equation}\label{za4}
{}_0^C D_x^{\varrho (x)} : = {}_0I_x^{n - q (x,r)}  \circ \frac{{d^n }}
{{dx^n }},
\end{equation}
where $q (x,r)=\varrho (x),\ q (x,r)=\varrho (r)$ and $q(x,r)=\varrho (x-r)$, in cases
\eqref{za1}-\eqref{za3}. Thus, we obtain, respectively:
\begin{enumerate}
  \item \textit{The type I: left Caputo fractional derivative of order $\varrho(x)$}
 \begin{equation}\label{2.2}
{}_0^{C} \mathscr{D}_x^{\varrho(x)}  u(x) = \frac{1}
{{\Gamma (n - \varrho (x))}}\int_0^x {\frac{{u^{(n)}(r)dr}}
{{(x - r)^{\varrho (x) - n + 1} }}}.
\end{equation}
  \item \textit{The type II: left Caputo fractional derivative of order $\varrho(x)$}
 \begin{equation}\label{2.2za}
{}_0^C \mathfrak{D}_x^{\varrho (x)} u(x) = \int_0^x {\frac{1}
{{\Gamma (n - \varrho (r))}}\frac{{u^{(n)} (r)dr}}
{{(x - s
)^{\varrho (r) - n + 1} }}}.
\end{equation}
   \item \textit{The type III: left Caputo fractional derivative of order $\varrho(x)$}
 \begin{equation}\label{2.2zaa}
{}_0^C \mathbb{D}_x^{\varrho (x)} u(x) = \int_0^x {\frac{1}
{{\Gamma (n - \varrho (x-r))}}\frac{{u^{(n)} (r)dr}}
{{(x - r)^{\varrho (x-r) - n + 1} }}},
\end{equation}
\end{enumerate}
 where $n-1< \varrho (\cdot)< n \in \mathbb{N}$. Such  operators have been used by researchers, for examples,
Coimbra et al. \cite{coimbra2003mechanics,soon2005variable} employed   the first type in the modeling of viscous-viscoelastic oscillator. Ingman and Suzdalnitsky \cite{ingman2005application} used the second type in the modeling of viscoelastic deformation process. Atanackovic and Pilipovic \cite{atanackovic2011hamilton} used  the third type in generalization of Hamilton's principle. Sun et al. \cite{sun2011comparative} introduced a comparative investigation of constant-order fractional derivative and the first two types of VO-F derivatives in characterizing the memory property of systems. However, the differences between the three types in applications are still not clear. There are other definitions of VO-F derivatives \cite{atanackovic2011hamilton}. In this paper, we will focus our attention on the first type of VO-F integrals and derivatives, i.e. $q(x,s)=\varrho(x)$.

The operator $ {}_0^{C} \mathscr{D}_x^{\varrho(x)}$ satisfies the following property $(n-1<  \varrho(x)\leq n \in \mathbb{N})$
 \begin{equation}\label{2.4} \begin{split}
&{}_0^{C} \mathscr{D}_x^{\varrho(x)} x^\gamma  = \begin{cases}
  \q \q 0, &  \gamma =0,\ldots, n-1,\\
   \dfrac {\Gamma{(\gamma +1)}}{\Gamma{(\gamma +1-\varrho(x))}}\ x^{\gamma -\varrho(x)}, \q\q\q & \gamma =n,n+1,\cdots.
  \end{cases}
\end{split} \end{equation}
Next, let us introduce some properties of the modified generalized Laguerre polynomials \cite{baleanu2013modified}.
Let $\Lambda =\{x \ | \ \ 0< x< \infty \}$ and $\chi(x)$  be a certain  weight  function on $\Lambda $ in the usual sense.
Define
\begin{equation}
  L_{\chi }^2 (\Lambda) =
  \{ y \ |\ y \ \textbf{is} \ \textbf{measurable} \ \& \ ||y||_{\chi }  < \infty \},
   \end{equation}
which is a Hilbert space, equipped with the following inner product and norm:
\begin{equation}(u,y)_{\chi }  = \int_{\Lambda } {\chi (x) u(x)y(x)}  \ dx,\q\q\ ||y||_{\chi }  = (y,y)_{\chi }^{\frac{1}{2}}. \end{equation}

In principle, the generalized Laguerre polynomials are suitable for the approximation of functions with growth at infinity.
We now recall some properties of the generalized Laguerre polynomials $\mathscr{L}_i^{(\theta ,\beta )} (x)$.

Let $\chi ^{(\theta ,\beta  )}  = x^\theta  e^{ - \beta x} ,\ \theta  >  - 1,$
 and $\beta>0$,  $\partial _x v(x) = \frac{\partial }
{{\partial x}}v(x)$. The corresponding generalized Laguerre polynomials of
degree $i$ are defined by
\begin{equation}\mathscr{L}_i^{(\theta ,\beta )} (x) = \frac{1}
{{i!}}x^{ - \theta } e^{\beta x} \partial _x^i (x^{i + \theta } e^{ - \beta x} ),\q\q i = 0,1,2, \ldots.\end{equation}
They are the eigenfunctions of the Sturm-Liouville problem
\begin{equation}\label{2.1}\partial _x (\chi ^{(\theta  + 1,\beta )} (x)\partial _x v(x)) + \lambda _i^{(\beta )} \chi ^{(\theta ,\beta )} (x)v(x) = 0,\q\ 0 < x < \infty,
\end{equation}
with the corresponding eigenvalues $\lambda _i^{(\beta )}=i \beta$. They fulfill the following three-term  recurrence relations:
\begin{equation}\label{2.2}\begin{split}
&\mathscr{L}_0^{(\theta ,\beta )} (x) = 1,\q\q \mathscr{L}_1^{(\theta ,\beta )} (x) =  - \beta x + \frac{{\Gamma (\theta  + 2)}}
{{\Gamma (\theta  + 1)}},\\
&\mathscr{L}_{i+1}^{(\theta ,\beta )} (x)  = \frac{{2i + \theta  + 1 - \beta x}}
{{i + 1}}\mathscr{L}_i^{(\theta ,\beta )} (x)  - \frac{{i + \theta }}
{{i + 1}}\mathscr{L}_{i-1}^{(\theta ,\beta )} (x),\q\q i \geq 1,
\end{split}\end{equation}
and
\begin{equation}\label{2.3}
\mathscr{L}_i^{(\theta ,\beta )} (x) =\mathscr{L}_i^{(\theta +1,\beta )} (x) -\mathscr{L}_{i-1}^{(\theta+1 ,\beta )} (x)  = \tfrac{1}
{\beta }(\partial _x \mathscr{L}_i^{(\theta ,\beta )} (x) - \partial _x \mathscr{L}_{i+1}^{(\theta ,\beta )} (x)  ).
\end{equation}
The  $m$-th derivative of a generalized Laguerre polynomial satisfies the relation
\begin{equation}\label{2.4}
\partial^m _x \mathscr{L}_i^{(\theta ,\beta )} (x)  = ( - \beta )^m \mathscr{L}_{i-m}^{(\theta+m ,\beta )} (x).
\end{equation}
The generalized Laguerre polynomials form a complete $L_{\chi^{\theta,\beta} }^2 (\Lambda)$-orthogonal system,
\begin{equation}\label{2.5}
\left( {\mathscr{L}_i^{(\theta ,\beta )} (x),\mathscr{L}_j^{(\theta ,\beta )} (x)} \right)_{\chi ^{(\theta ,\beta )} }  = \gamma_i ^{(\theta ,\beta )} \delta _{i,j},\q\q \gamma_i ^{(\theta ,\beta )}  = \frac{{\Gamma (i + \theta  + 1)}}
{{\beta ^{\theta  + 1} \Gamma (i + 1)}},
\end{equation}
where $\delta _{i,j}$ is the Kronecker symbol.
\section{Approximation to the variable-order fractional integral} \label{sec3}
The main goal of this section is to  develop a new algorithm to approximate the VO-F integral of a given function.

Let $u \in L_{\chi^{\theta,\beta} }^2(\Lambda)$ and $N$ be a
positive integer, then, we can expand it in terms of generalized Laguerre polynomials as
\begin{equation}\label{2.6}
u(x) \simeq u_N (x) = \sum\limits_{i = 0}^N {\ell _i^{(\theta ,\beta )} \mathscr{L}_i^{(\theta ,\beta )} (x)}.
\end{equation}
If $u_N (x)$ is an orthogonal projection of $u(x)$, then $\ell _i^{(\theta ,\beta )}$ can be
determined by the orthogonality of $\{ \mathscr{L}_i^{(\theta ,\beta )} (x) \}$. If $u_N (x)$ is the interpolation of $u(x)$ on
the generalized Laguerre-Gauss points $\{ x^{\theta,\beta}_{N,i} \}_{i=0}^N$, then $\ell _i^{(\theta ,\beta )}$ can be
determined by
\begin{equation}\label{2.7}
\ell _i^{(\theta ,\beta )} = \frac{1}
{{\gamma _i^{(\theta ,\beta )} }}\sum\limits_{j = 0}^N {u(x_{N,j}^{(\theta ,\beta )} )} \mathscr{L}_i^{(\theta ,\beta )} (x_{N,j}^{(\theta ,\beta )} )\chi ^{(\theta ,\beta )}_{N,j},
\end{equation}
 where
 $\chi ^{(\theta ,\beta )}_{N,j}$ are the corresponding weights. In this paper, we assume that $u_N (x)$ is the interpolation of $u(x)$.

Therefore, for any $n-1<\varrho_{min}<\varrho(x)<\varrho_{max}<n$, the VO-F integral $ {}_0I_x^{\varrho (x)} u(x)$ can be approximated by
   \begin{equation}\label{2.9}\begin{split}
 {}_0I_x^{\varrho (x)} u(x)& \approx  {}_0I_x^{\varrho (x)} u_N (x) = \frac{1}
{{\Gamma (\varrho (x))}}\int_0^x {(x - t)^{\varrho (x) - 1} u_N (t)dt}\\
&=\frac{1}
{{\Gamma (\varrho (x))}}\sum\limits_{i = 0}^N {\ell _i^{(\theta ,\beta )}  \int_0^x {(x - t)^{\varrho (x) - 1} \mathscr{L}_i^{(\theta ,\beta )} (t)dt} }  \\
&= \sum\limits_{i = 0}^N {\ell _i^{(\theta ,\beta )}  \mathscr{\hat{L}}_i^{(\varrho(x),\theta ,\beta )} (x)},
 \end{split}\end{equation}
where
\begin{equation} \mathscr{\hat{L}}_i^{(\varrho(x),\theta ,\beta )} (x)= \frac{1}
{{\Gamma (\varrho (x))}}\int_0^x {(x - t)^{\varrho (x) - 1} \mathscr{L}_i^{(\theta ,\beta )} (t)dt}\end{equation}.

 Next, we develop a recurrence formula to calculate the $\varrho(x)$th-order fractional integral of $ \mathscr{L}_i^{(\theta ,\beta )}$.

From \eqref{2.2}, we can easily get
\begin{equation}
\mathscr{\hat{L}}_0^{(\varrho(x),\theta ,\beta )} (x)=\frac{{x^{\varrho (x)} }}
{{\Gamma (\varrho (x) + 1)}},
\end{equation}
and
\begin{equation}
 \mathscr{\hat{L}}_1^{(\varrho(x),\theta ,\beta )} (x)=\frac{{x^{\varrho (x)} }}
{{\Gamma (\varrho (x) + 1)}}-\frac{{\beta x^{\varrho (x)+1} }}
{{\Gamma (\varrho (x) + 2)}}.
\end{equation}
For $i\geq1$, we have
 \begin{equation}\label{2.10}\begin{split}
\mathscr{\hat{L}}_{i+1}^{(\varrho(x),\theta ,\beta )} (x)&=\frac{1}{{\Gamma (\varrho (x))}}\int_0^x {(x - t)^{\varrho (x) - 1}} \mathscr{L}_{i+1}^{(\theta ,\beta )} (t)dt\\
&=\frac{1}{{\Gamma (\varrho (x))}}\frac{1}{{i + 1}}\int_0^x {(x - t)^{\varrho (x) - 1}} \\
&\times \q\left[ (2i + \theta  + 1 - \beta t)\mathscr{L}_i^{(\theta ,\beta )} (t)  - (i + \theta)\mathscr{L}_{i-1}^{(\theta ,\beta )} (t)\right]dt\\
&=\frac{1}{{i + 1}} \Big\{(2i + \theta  + 1 - \beta x)\mathscr{\hat{L}}_i^{(\varrho (x),\theta ,\beta )} (x)  - (i + \theta)\mathscr{\hat{L}}_{i-1}^{(\varrho (x),\theta ,\beta )} (x)\\
&+\frac{\beta}{{\Gamma (\varrho (x))}}\int_0^x {(x - t)^{\varrho (x) }} \mathscr{L}_{i}^{(\theta ,\beta )} (t)dt \Big\}.
 \end{split}\end{equation}
 We can verify from formula \eqref{2.3} that for $i\geq1$,
 \begin{equation}\label{2.11}\begin{split}
\mathscr{\hat{L}}_{i+1}^{(\varrho(x),\theta ,\beta )} (x)=&\frac{1}{{i + 1}} \Big\{(2i + \theta  + 1 - \beta x)\mathscr{\hat{L}}_i^{(\varrho (x),\theta ,\beta )} (x)  - (i + \theta)\mathscr{\hat{L}}_{i-1}^{(\varrho (x),\theta ,\beta )} (x)\\
&+\frac{1}{{\Gamma (\varrho (x))}}\int_0^x {(x - t)^{\varrho (x) }} \big(\partial _t \mathscr{L}_i^{(\theta ,\beta )} (t) - \partial _t \mathscr{L}_{i+1}^{(\theta ,\beta )} (t)  \big)dt \Big\}\\
=&\frac{1}{{i + 1}} \Big\{(2i + \theta  + 1 - \beta x)\mathscr{\hat{L}}_i^{(\varrho (x),\theta ,\beta )} (x)  - (i + \theta)\mathscr{\hat{L}}_{i-1}^{(\varrho (x),\theta ,\beta )} (x)\\
&+\frac{1}{{\Gamma (\varrho (x))}}\Big[{(x - t)^{\varrho (x) }} \big( \mathscr{L}_i^{(\theta ,\beta )} (t) - \mathscr{L}_{i+1}^{(\theta ,\beta )} (t) \big)  \Big]_0^x \\
&+\varrho(x) \big( \mathscr{\hat{L}}_i^{(\varrho(x),\theta ,\beta )} (x) - \mathscr{\hat{L}}_{i+1}^{(\varrho(x),\theta ,\beta )} (x) \big)\Big\}\\
=&\frac{1}{{(i + \varrho(x)+1)}} \Big\{(2i + \theta  +\varrho (x)+ 1 - \beta x)\mathscr{\hat{L}}_i^{(\varrho (x),\theta ,\beta )} (x) \\
& - (i + \theta)\mathscr{\hat{L}}_{i-1}^{(\varrho (x),\theta ,\beta )} (x)-\frac{x^{\varrho(x)}}{{\Gamma (\varrho (x))}} \big( \mathscr{L}_i^{(\theta ,\beta )} (0) - \mathscr{L}_{i+1}^{(\theta ,\beta )} (0) \big)   \Big\}.
 \end{split}\end{equation}
Hence, for $i\geq1$, we get the following recurrence relation
  \begin{equation}\label{2.12}\begin{split}
\mathscr{\hat{L}}_{i+1}^{(\varrho(x),\theta ,\beta )} (x)=&\frac{1}{{(i + \varrho(x)+1)}} \Big\{(2i + \theta  +\varrho (x)+ 1 - \beta x)\mathscr{\hat{L}}_i^{(\varrho (x),\theta ,\beta )} (x) \\
& - (i + \theta)\mathscr{\hat{L}}_{i-1}^{(\varrho (x),\theta ,\beta )} (x)-\frac{x^{\varrho(x)}}{{\Gamma (\varrho (x))}} \big( \mathscr{L}_i^{(\theta ,\beta )} (0) - \mathscr{L}_{i+1}^{(\theta ,\beta )} (0) \big)   \Big\}.
 \end{split}\end{equation}
 So, $\mathscr{\hat{L}}_{i}^{(\varrho(x),\theta ,\beta )} (x)$ can be calculated by the following formula
  \begin{equation}\label{2.13}
\left\{  \begin{split}
  &\mathscr{\hat{L}}_0^{(\varrho(x),\theta ,\beta )} (x)=\frac{{x^{\varrho (x)} }}
{{\Gamma (\varrho (x) + 1)}},  \\
 & \mathscr{\hat{L}}_1^{(\varrho(x),\theta ,\beta )} (x)=\frac{{x^{\varrho (x)} }}
{{\Gamma (\varrho (x) + 1)}}-\frac{{\beta x^{\varrho (x)+1} }}
{{\Gamma (\varrho (x) + 2)}},  \\
 & \begin{split} \mathscr{\hat{L}}_{i+1}^{(\varrho(x),\theta ,\beta )} (x)=&\frac{1}{{(i + \varrho(x)+1)}} \Big\{(2i + \theta  +\varrho (x)+ 1 - \beta x)\mathscr{\hat{L}}_i^{(\varrho (x),\theta ,\beta )} (x) \\
  &- (i + \theta)\mathscr{\hat{L}}_{i-1}^{(\varrho (x),\theta ,\beta )} (x)-\frac{x^{\varrho(x)}}{{\Gamma (\varrho (x))}} \big( \mathscr{L}_i^{(\theta ,\beta )} (0) - \mathscr{L}_{i+1}^{(\theta ,\beta )} (0) \big)   \Big\},  \\
 &\qq\qq  i\geq1. \end{split}
\end{split} \right.\end{equation}
Therefore, $ {}_0I_x^{\varrho (x)}u(x)$ can be approximated  by
 \begin{equation}\label{2.14}
 {}_0I_x^{\varrho (x)} u(x) \approx  {}_0I_x^{\varrho (x)} u_N (x)   = \sum\limits_{i = 0}^N {\ell _i^{(\theta ,\beta )}  \mathscr{\hat{L}}_i^{(\varrho(x),\theta ,\beta )} (x)},
\end{equation}
 where $\mathscr{\hat{L}}_i^{(\varrho(x),\theta ,\beta )} (x)$ is given by \eqref{2.13}, and  $\ell _i^{(\theta ,\beta )}$ is given by \eqref{2.7}.
 \begin{rem}
 When $\varrho(x)=\varrho=constant$, then the operator of
VO-F integral is reduced to corresponding
integral of constant order and the approximation relation \eqref{2.14} reduces to
\begin{equation*}
 {}_0I_x^{\varrho } u(x) \approx  {}_0I_x^{\varrho } u_N (x)   = \sum\limits_{i = 0}^N {\ell _i^{(\theta ,\beta )}  \mathscr{\hat{L}}_i^{(\varrho,\theta ,\beta )} (x)}.
\end{equation*}
\end{rem}
\section{Approximation to the   variable-order fractional Caputo derivative} \label{sec4}
In this section, we describe how to use the  generalized Laguerre polynomials to numerically approximate the VO-F derivative in the Caputo sense of a given function $u(x), x \in \Lambda$. The computerized mathematical algorithm is based on the
numerical approximation of the VO-F integral derived in the previous
section.

Suppose that $u_N(x)$ is the approximate polynomial of $u(x)$, which can be expressed by
\begin{equation}\label{4.1}
 u_N (x) = \sum\limits_{i = 0}^N {\ell _i^{(\theta ,\beta )} \mathscr{L}_i^{(\theta ,\beta )} (x)},\q x \in \Lambda.
\end{equation}
Let $n-1<\varrho_{min}<\varrho(x)<\varrho_{max}<n\in \mathbb{N}$, we approximate the VO-F derivative ${}_0^{C}\mathscr{D}_x^{  \varrho (x)} u(x)$ by the generalized Laguerre polynomials as
   \begin{equation}\label{4.2}\begin{split}
{}_0^C\mathscr{D}_x^{ \varrho (x)} u_N (x) &= \frac{1}
{{\Gamma (n-\varrho (x))}}\int_0^x {(x - t)^{n-\varrho (x) - 1}\partial^n_t u_N (t)dt}\\
&=\frac{1}
{{\Gamma (n-\varrho (x))}}\sum\limits_{i = 0}^N {\ell _i^{(\theta ,\beta )}  \int_0^x {(x - t)^{n-\varrho (x) - 1} \partial^n_t\mathscr{L}_i^{(\theta ,\beta )} (t)dt} }\\
&=\sum\limits_{i = 0}^N {\ell _i^{(\theta ,\beta )}(-\beta)^n \Big( \frac{1}
{{\Gamma (n-\varrho (x))}}\int_0^x {(x - t)^{n-\varrho (x) - 1} \mathscr{L}_{i-n}^{(\theta+n ,\beta )} (t)dt} \Big)} \\
&= \sum\limits_{i = 0}^N {\ell _i^{(\theta ,\beta )}(-\beta)^n   \mathscr{\hat{L}}_{i-n}^{(n-\varrho(x),\theta+n ,\beta )} (x)}.
 \end{split}\end{equation}
where $\mathscr{\hat{L}}_{i-n}^{(n-\varrho(x),\theta+n ,\beta )} (x)=0$ for $0\leq i \leq n-1$.

Therefore, the VO-F derivative of $u(x)$ can be approximated by
\begin{equation}\label{4.3}
{}_0^C\mathscr{D}_x^{  \varrho (x)} u_N (x)   = \sum\limits_{i = 0}^N {\ell _i^{(\theta ,\beta )}(-\beta)^n   \mathscr{\hat{L}}_{i-n}^{(n-\varrho(x),\theta+n ,\beta )} (x)}=\sum\limits_{i = 0}^N {\ell _i^{(\theta ,\beta )}D_{i,n,\theta ,\beta }^{(\varrho (x))}(x)},
\end{equation}
which alternatively may be written in the matrix form:
\begin{equation}\label{4.4}
{}_0^C\mathscr{D}_x^{  \varrho (x)} u_N (x)   =\mathcal{L} ^{(\theta ,\beta )}\textbf{D}_{n,\theta ,\beta }^{(\varrho (x))}(x),
\end{equation}
with
 \begin{equation}\label{4.5}\begin{split}
 \mathcal{L} ^{(\theta ,\beta )}& = \left[ {\ell _0^{(\theta ,\beta )},\ell _1^{(\theta ,\beta )}, \ldots ,\ell _N^{(\theta ,\beta )}} \right],\\
\textbf{D}_{n,\theta ,\beta }^{(\varrho (x))}(x)&= \left[ {D_{0,n,\theta ,\beta }^{(\varrho (x))}(x),D_{1,n,\theta ,\beta }^{(\varrho (x))}(x), \ldots ,D_{N,n,\theta ,\beta }^{(\varrho (x))}(x)} \right]^T,
 \end{split}\end{equation}
where $n-1<\varrho_{min}<\varrho(x)<\varrho_{max}<n\in \mathbb{N},$  and
 \begin{equation}\label{4.6}\begin{split}
D_{i,n,\theta ,\beta }^{(\varrho (x))}(x)&=\frac{1}
{{\Gamma (n-\varrho (x))}}\int_0^x {(x - t)^{n-\varrho (x) - 1} \partial^n_t\mathscr{L}_i^{(\theta ,\beta )} (t)dt}\\
&=(-\beta)^n   \mathscr{\hat{L}}_{i-n}^{(n-\varrho(x),\theta+n ,\beta )}, \q i\geq n,
 \end{split}\end{equation}
and $D_{i,n,\theta ,\beta }^{(\varrho (x))}(x)=0$ for $0\leq i \leq n-1$.
\section{Applications of the algorithms} \label{sec5}
After the construction of the VO-F differentiation matrix of Caputo type  \eqref{4.4}, we now use
the generalized Laguerre spectral collocation method together with the generalized Laguerre operational
matrix of VO-F derivative to solve the following VO-F differential equation:
 \begin{equation}\label{5.1}\begin{split}
&a(x)u'(x) + b(x){}_0^{C}\mathscr{D}_x^{  \varrho (x)}u(x) + c(x)u(x) = f(x),\q\q x \in \Lambda,\\
&u(0)=u_0.
 \end{split}\end{equation}
where $0<\varrho_{min}<\varrho(x)<\varrho_{max}<1,\ a(x),\ b(x),\ c(x)$ and $f(x)$ are real-valued functions. A special case of \eqref{5.1} occurs when $\varrho(x)=\varrho=\frac{1}{2}$, corresponding to the fractional Basset equation. This model represents a classical problem in fluid dynamics in the
scope of an unsteady motion of a particle that accelerates in a viscous fluid under the action of gravity.

Suppose that $u_N(x)$ is the approximate solution of $u(x)$ and $x_i\ (0\leq i \leq N)$ is the generalized Laguerre-Gauss nodes of $\mathscr{L}_{N+1}^{(\theta ,\beta )} (x)$
Now, using \eqref{4.4} then it is easy to write
\begin{equation}\label{5.2}\begin{split}
&a(x_i)\mathcal{L} ^{(\theta ,\beta )}\textbf{D}_{1,\theta ,\beta }^{(1)}(x_i)
 + b(x_i)\mathcal{L} ^{(\theta ,\beta )}\textbf{D}_{1,\theta ,\beta }^{(\varrho (x_i))}(x_i)
 + c(x_i)\mathcal{L} ^{(\theta ,\beta )}\textbf{D}_{0,\theta ,\beta }^{(0)}(x_i)
 = f(x_i),\\
&u_N(0)=u_0 \q\q \q i=0,1,\ldots,N-1.
 \end{split}\end{equation}
Let us denote $F=\left[ {f(x_0),\ldots,f(x_{N-1}),u_0} \right]$. Then \eqref{5.2}  is equivalent to the following matrix equation
\begin{equation}\mathcal{L} ^{(\theta ,\beta )} E = F,\end{equation}
where
\begin{equation}\label{5.3}
E=(e_{ij})  = \left\{   \begin{split}
 & \begin{split} &a(x_j)D_{i,1,\theta ,\beta }^{(1)}(x_j)
 + b(x_j)D_{i,1,\theta ,\beta }^{(\varrho (x_j))}(x_j),
 + c(x_j)D_{i,0,\theta ,\beta }^{(0)}(x_j) \\
 & \q\q \q 0\leq i \leq   N,\ 0\leq j \leq N-1, \end{split} \\
 & \frac{{\Gamma (i + \theta  + 1)}}
{{\Gamma (1 + \theta )i!}},\q\q\ 0 \leq i \leq N,\ j = N.
\end{split}  \right.
\end{equation}

For the following fractional initial value problem
\begin{equation}\label{5.4}\begin{split}
&a(x)u^{(m)}(x) + b(x){}_0^{C}\mathscr{D}_x^{  \varrho (x)}u(x) + c(x)u(x) = f(x),\q\q x \in \Lambda,\\
&u(0)=u_0,\q\q\q u'(0)= v_0,
 \end{split}\end{equation}
 where $1<\varrho_{min}<\varrho(x)<\varrho_{max}<2,\ m=1$ or $2$, we can also get the algebraic equation of
the form like \eqref{5.2}, where $F=\left[ {f(x_0),\ldots,f(x_{N-2}),u_0,v_0} \right]$ and
\begin{equation}\label{5.5}
E=(e_{ij})  = \left\{   \begin{split}
 & \begin{split} &a(x_j)D_{i,m,\theta ,\beta }^{(m)}(x_j)
 + b(x_j)D_{i,2,\theta ,\beta }^{(\varrho (x_j))}(x_j),
 + c(x_j)D_{i,0,\theta ,\beta }^{(0)}(x_j) \\
 & \q\q \q 0\leq i \leq   N,\ 0\leq j \leq N-2, \end{split} \\
 & \frac{{\Gamma (i + \theta  + 1)}}
{{\Gamma (1 + \theta )i!}},\q\q\ 0 \leq i \leq N,\ j = N-1,\\
&D_{i,1,\theta ,\beta }^{(1)}(0),\q\q\ 0 \leq i \leq N,\ j = N.
\end{split}  \right.
\end{equation}
In our implementation, these systems have been solved using the Mathematica function FindRoot with zero initial approximation.

A special case of \eqref{5.4} occurs when $m=2$ and $\varrho(x)=\varrho=\frac{3}{2}$, corresponding to the fractional Bagley-Torvik equation. The fractional Bagley-Torvik equation is important for modeling
the motion of a thin rigid plate immersed in a Newtonian fluid.

  \section{Numerical examples} \label{sec6}
\begin{figure}[h!]\centering
\includegraphics[width=0.8\textwidth]{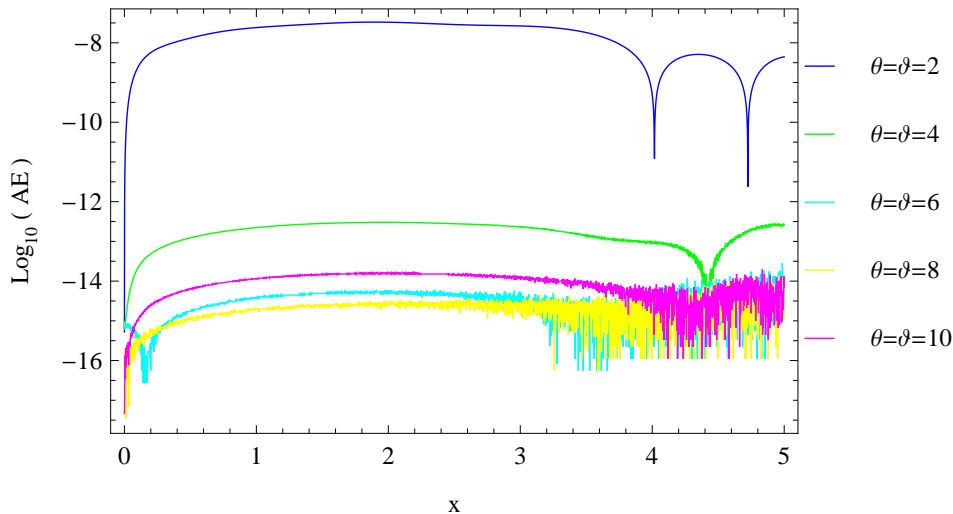}
\caption{Logarithmic graphs of the AE at $\varrho(x)=\frac{{9 + \sin (x - 10)}}{5}$, $N=25$ and various choices of $\theta=\beta$ for Example \ref{Ex.2}.}
\label{fig3}\end{figure}

  In this section we  present three numerical examples to illustrate the  accuracy and
efficiency of the
algorithms presented in the previous sections. The first one is introduced to test the accuracy of the formula \eqref{4.4}.
\begin{example}\label{Ex.1}
Let $u(x)=e^x,\ x \in [0,1]$. Now we numerically calculate the VO-F derivative ${}_0^{C}\mathscr{D}_x^{\varrho (x)}u(x),\ n-1 \leq \varrho(x) \leq n.$

The analytical form of the VO-F derivative of $u(x)$ is given by
\[{}_0^{C}\mathscr{D}_x^{\varrho (x)}e^x=e^x \left(1 - \frac{{\Gamma (n - \varrho (x),x)}}{{\Gamma (n - \varrho (x))}} \right).\]
\end{example}
In Table \ref{tab5.1} and \ref{tab5.2}, we list the maximum absolute errors (AEs) between the exact solution $u(x)$ and the approximate
solutions $u_N$ at different constant- and variable-orders, respectively.
\begin{table} [h!] \centering
  \caption{The maximum AEs  at different constant fractional orders for Example \ref{Ex.1}. }\label{tab5.1}
  \centerline{ \begin{tabular}{cccccccccc}
                   \hline
                     $N$ & $(\theta,\beta)$& $\varrho=0.2$   &  $\varrho=0.5$  & $\varrho=0.8$  &  $\varrho=1.2$ &  $\varrho=1.5$   &  $\varrho=1.8$  \\
                    \hline
                  10& (1,3)&$7.93\times10^{-3} $   &$1.46\times10^{-2} $&$3.36\times10^{-2}$&$1.10\times10^{-1}$&$1.78\times10^{-1}$&$3.62\times10^{-1}$\\
                  20&      &$1.53\times10^{-5} $   &$3.45\times10^{-5} $&$9.72\times10^{-5}$&$4.10\times10^{-4}$&$8.24\times10^{-4}$&$2.07\times10^{-3}$\\
                  40&      &$3.07\times10^{-11} $  &$8.49\times10^{-11} $&$2.93\times10^{-10}$&$1.61\times10^{-9}$&$4.00\times10^{-9}$&$1.25\times10^{-8}$\\
                  80&      &$4.88\times10^{-15} $  &$6.21\times10^{-15}$&$5.77\times10^{-15}$&$7.32\times10^{-15}$&$1.55\times10^{-14}$&$1.37\times10^{-14}$\\
                  \hline
                  10& (2,6)&$2.93\times10^{-6} $  &$6.04\times10^{-6}$&$1.55\times10^{-5}$&$5.76\times10^{-5}$&$1.06\times10^{-4}$&$2.48\times10^{-4}$\\
                  20&   &$1.09\times10^{-12} $  &$2.73\times10^{-12}$&$8.64\times10^{-12}$&$4.13\times10^{-11}$&$9.46\times10^{-11}$&$2.73\times10^{-10}$\\
                    40&    &$1.33\times10^{-15} $ &$2.66\times10^{-15}$&$2.67\times10^{-15}$&$1.77\times10^{-15}$&$3.10\times10^{-15}$&$2.66\times10^{-15}$\\
                     80&    &$1.33\times10^{-15}$&$2.66\times10^{-15}$&$2.67\times10^{-15}$&$1.77\times10^{-15}$&$3.10\times10^{-15}$&$2.66\times10^{-15}$\\
                    \hline
                      \end{tabular}}
                     \end{table}
                     \begin{table}[h!]  \centering
  \caption{The maximum AEs at different variable orders for Example \ref{Ex.1}. }\label{tab5.2}
  \centerline{ \begin{tabular}{cccc}
                 \hline
               $(\theta,\beta) $& $N$ & $\varrho(x)=\frac{{9 + \sin t}}{{10}}$ & $\varrho(x)=\frac{{3 + \tanh t}}{2}$ \\
               \hline
        (2,4)   & $10$ & $4.648\times10^{-3}$ & $1.833\times10^{-2}$ \\
                 & 20 & $4.556\times10^{-7}$  &  $2.598\times10^{-6}$  \\
                 & 30 & $2.282\times10^{-11}$  & $1.625\times10^{-10}$ \\
                 & 40 & $5.329\times10^{-15}$  & $7.688\times10^{-15}$ \\
                 \hline
                 (3,6) & 10 & $9.862\times10^{-5}$ & $4.228\times10^{-4}$ \\
                  & 20 &$1.013\times10^{-10}$  & $6.287\times10^{-10}$  \\
                   & 30&$3.997\times10^{-15}$  &$3.552\times10^{-15}$  \\
                  & 40 &$3.997\times10^{-15}$  & $3.552\times10^{-15}$  \\
                 \hline
               \end{tabular}
  }
                     \end{table}
\begin{example}\label{Ex.2} Consider the following VO-F Bagley-Torvik equation
                  \begin{equation}\label{2.017}\begin{split}
                  &u''(x) + {}_0^{C}\mathscr{D}_x^{  \varrho (x)}u(x) + u(x) = f(x),\q\q x \in (0,L],\\
                  &u(0)=0\q\q u'(0)=1.
                   \end{split}\end{equation}
The exact solution is $u(x)=\sin x$.
   \end{example}
   Table \ref{tab5.3} displays the maximum AEs of the proposed method at $\varrho(x)=\frac{3}{2}$ and $\varrho(x)=\frac{{9 + \sin (x - 10)}}{5}$ for different chooses of $\theta,\ \beta$ and $N$. Figure \ref{fig3} shows the logarithmic graphs of the AEs  at $\varrho(x)=\frac{{9 + \sin (x - 10)}}{5}$, $N=25$ and different chooses of $\theta=\beta$.
\begin{table} [h!] \centering
  \caption{The maximum AEs at $L=1$ for Example \ref{Ex.2}.}\label{tab5.3}
  \centerline{ \begin{tabular}{ccccc}
                 \hline
               $\varrho(x)  $& $N$ &$\theta=0,\beta=1$ &$\theta=2,\beta=4$ & $\theta=3,\beta=6$ \\
               \hline
       $\frac{3}{2}$    & $5$ & $5.546\times10^{-3}$  &$2.916\times10^{-4}$ &$1.427\times10^{-4}$ \\
                       & 10  & $4.485\times10^{-4}$   &$1.431\times10^{-7}$ &$9.038\times10^{-9}$   \\
                         & 15  & $8.845\times10^{-6}$&$3.675\times10^{-11}$&$9.313\times10^{-12}$ \\
                        & 20  & $8.133\times10^{-6}$  &$2.166\times10^{-13}$&$2.220\times10^{-15}$ \\
                 \hline
  $\frac{{9 + \sin (x - 10)}}{5}$ & 5 &$8.318\times10^{-3}$ & $2.666\times10^{-3}$ & $1.231\times10^{-3}$ \\
                                 & 10 &$2.515\times10^{-3}$ &$3.854\times10^{-6}$  & $1.179\times10^{-7}$  \\
                                  & 15&$1.771\times10^{-4}$&$2.721\times10^{-9}$  &$7.242\times10^{-11}$  \\
                                 & 20 &$5.418\times10^{-6}$&$1.0522\times10^{-12}$  & $2.742\times10^{-14}$  \\
                 \hline
               \end{tabular}
  }
                     \end{table}
                   \begin{example}\label{Ex.3} Consider the following VO-F Bagley-Torvik equation
                  \begin{equation}\label{2.017}\begin{split}
                  &u''(x) + {}_0^{C}\mathscr{D}_x^{  \varrho (x)}u(x) + u(x) = \frac{{\Gamma (4)}}{{\Gamma (4 - \varrho (x))}}x^{3 - \varrho (x)}  + x^3  + 7x + 1,\ x\in (0,\frac{\pi}{2}],\\
                  &u(0)= u'(0)=1.
                   \end{split}\end{equation}
The exact solution is $u(x)=x^3  + x + 1$.
   \end{example}
   In this example, we consider two cases, $\varrho(x) = 1.5$ and $\varrho(x) =1 + 0.5\left| {\sin x} \right|$, $x \in [0,\frac{\pi}{2}]$. Table \ref{t200} displays a comparison between
the M- and IM-algorithms \cite{moghaddam2017integro,moghaddam2017extended} by means of the maximum AEs.
\begin{table}[h!]
\caption{Comparison of the maximum AEs for  Example \ref{Ex.3} with $\theta=\beta=10$ and two cases of $\varrho(x)$, $x \in [0,\frac{\pi}{2}]$. }\label{t200}
                  \begin{tabular}{c|ccc|ccccccccccc}
                   \hline
                   $\varrho(x)$& &M-algorithm \cite{moghaddam2017extended} &IM-algorithm\cite{moghaddam2017integro} &&Collocation method   \\
                    \hline
                                     &$h=0.02$ &$8.16\times10^{-3}$ &$1.15\times10^{-3}$&$N=3$&$5.77\times10^{-15}$&\\
                    $1.5$            &$h=0.01$ &$4.27\times10^{-3}$ &$4.29\times10^{-4}$&$N=4$&$4.57\times10^{-15}$&    \\
                                     &$h=0.005$&$2.18\times10^{-3}$ &$1.56\times10^{-4}$&$N=5$&$4.44\times10^{-15}$&   \\
                    \hline
                                     &$h=0.02$ &$7.22\times10^{-3}$&$ 1.53\times10^{-4}$&$N=3$&$4.88\times10^{-15}$&\\
    $1 + 0.5\left| {\sin x} \right|$ &$h=0.01$ &$3.70\times10^{-3}$&$ 5.57\times10^{-5}$&$N=4$&$3.10\times10^{-15}$&    \\
                                     &$h=0.005$&$1.86\times10^{-3}$&$ 1.93\times10^{-5}$&$N=5$&$2.77\times10^{-15}$&   \\
                   \hline
                  \end{tabular}
\end{table}
  \section{Conclusion}\label{sec7}

  In this paper, an efficient three-term recurrence relation to calculate both VO-F  integrals and derivatives of the modified generalized Laguerre polynomials was developed. Spectral collocation methods were developed to solve VO-F differential equations. The
results of this paper expand the application of the  Laguerre-Gauss collocation methods to VO-F problems. The suggested algorithms can be used also for solving VO-F fractional partial differential equations. Hence, the method is promising for VO-F differential equations. However, the analysis
of the scheme is a challenging problem deserving further study.


\end{document}